\documentclass[12pt]{amsart}

\usepackage{amsmath, amssymb, amsthm}
\usepackage{mathrsfs}

\usepackage{graphicx}
\usepackage{hyperref} 

\theoremstyle{plain}
\newtheorem{theorem}{Theorem}[section]
\newtheorem{lemma}[theorem]{Lemma}

\theoremstyle{definition}
\newtheorem{definition}[theorem]{Definition}
\newtheorem{example}[theorem]{Example}

\theoremstyle{remark}
\newtheorem{remark}[theorem]{Remark}

\newcommand{\R}{{\mathbb R}}
\newcommand{\C}{{\mathbb C}}

\newcommand{\N}{{\mathbb N}}
\newcommand{\Z}{{\mathbb Z}}
\newcommand{\D}{\mathbb{D}}
\newcommand{\ds}{\displaystyle}
\newcommand{\cal}{\mathcal}
\renewcommand{\Re}{\operatorname{Re}}
\renewcommand{\Im}{\operatorname{Im}}

\numberwithin{equation}{section}

\title[Boundary limits of the Bergman kernel and metric]
{Boundary limits of the Bergman kernel and metric on generally decoupled domains of infinite type}

\author{Rimu Shimogai}
\address{Faculty of Mathematics, Kyushu University, 
Motooka 744, Nishi-ku, Fukuoka, 819-0395, Japan}
\email{shimogai.rimu.286@s.kyushu-u.ac.jp}


\date{\today}

\begin{document}

\begin{abstract}
In this paper, we compute the precise nontangential boundary limits of the Bergman kernel on the diagonal and its related invariants on generally decoupled domains in $\C^{n+1}$  ($n\ge1$) at a boundary point of infinite type. Our results extend recent work by Jaiswal for exponentially flat domains.
\end{abstract}

\maketitle

\section{Introduction}
Let $D$ be a bounded domain in $\C^{n+1}$ ($n\ge1$) with smooth boundary $bD$, and let $\cal K_D$ be the Bergman kernel of $D$ (on the diagonal). The boundary behavior of $\cal K_D$ is one of the central topics in the study of complex analysis and geometry. When $D$ is strongly pseudoconvex, H\"ormander \cite{Hor65} proved that 
\begin{equation*}
\ds\lim_{z\to p\in bD}\cal K_D(z)d(z)^{n+2}=\ds\frac{(n+1)!}{4\pi^{n+1}}\cal L(p),
\end{equation*}
where $d(z)$ is the Euclidean distance from $z\in D$ to $bD$ and $\cal L(p)$ is the determinant of the Levi form at $p$. In this setting, Diederich \cite{Die70, Die73} computed the boundary limit of the Bergman metric, and later Klembeck \cite{Klem78} provided that of its holomorphic sectional curvature by using the asymptotic expansion of the Bergman kernel given by Fefferman \cite{Feff78}. Although, in this paper, we focus on the boundary limits in the sense of H\"ormander-type result of the Bergman invariants on weakly pseudoconvex domains, there have been many interesting results about lower and upper bounds of these invariants near the boundary points (see, for example, the references in \cite{BSY95} for the case of finite type, and \cite{Bha10, Bha20, Jais255} for the case of infinite type).

In the weakly pseudoconvex case, results, analogous to those in the strongly pseudoconvex case, have been obtained at least for certain special classes of domains. Boas, Straube and Yu \cite{BSY95} gave the boundary limits of the Bergman kernel and its associated invariants for a large class of domains of finite type, called {\it h}-extendible domains. Kamimoto \cite{Kami04, Kami24} investigated the asymptotic behavior of the Bergman kernel at the boundary for certain pseudoconvex model domains of not only finite but also infinite type, by using the geometrical information of the Newton polyhedron (a tool from singularity theory) determined by the defining function.

In the infinite type case, Jaiswal \cite{Jais254} computed the boundary limits of the Bergman kernel and its associated invariants on exponentially flat domains whose boundary contains the origin in $\C^{n+1}$ (see also \cite{KL02} for the case of dimension two). This domain is of infinite type that generalizes a domain whose local defining function near the origin is of the form $e^{-1/|z'|^2}-\Re z_1$, $(z_1,z')\in\C\times\C^n$.

In this paper, we focus on domains of infinite type with a certain decoupled structure in a neighborhood of the origin in $\C^{n+1}$ on the boundary, which are defined as follows:

\begin{definition}\label{D13}
Let $D\subset\C^{n+1}$ be a smoothly bounded domain with $0\in bD$. We say that $D$ is {\it generally decoupled near $0$} if there exists a small open neighborhood $U\subset\C^{n+1}$ of $0$ such that
\begin{equation}
D\cap U=\left\{(z_1,z')\in U : \rho(z_1,z')<0\right\},
\end{equation}
where $\rho$ is a local defining function of $D$ near $0$ taking the form
\begin{equation}\label{eq12}
\rho(z_1,z')=\ds\sum_{i=1}^kf_i(|z^i|)-\Re z_1,
\end{equation}
$z'=(z^1,\dots,z^k)\in\C^{n_1}\times\dots\times\C^{n_k}\ (k\ge1), n_1,\dots,n_k\in\N$ with $n_1+\dots+n_k=n$ and each of $f_i$ is either finite type at $0$ or exponentially flat at $0$. For the definition of these functions, see Section 2.1.
\end{definition}

This domain is a generalization of the {\it decoupled domains} considered by McNeal \cite{Mc91}, where he studied lower and upper bounds of the Bergman kernel and metric near the boundary points of finite type of those domains in $\C^n$. Furthermore, Jaiswal investigated the optimal bounds of the Bergman kernel and metric on generally decoupled domains (referred to as {\it generalized decoupled domains} in \cite{Jais255}). However, to the best of my knowledge, the precise boundary limits of these invariants of their domains have not yet been obtained. The goal of the present paper is to compute the boundary limits of the Bergman kernel, the Bergman metric, and its holomorphic sectional curvature on generally decoupled domains.

Before stating our main results, we introduce some notations. A {\it nontangential cone stream} approaching $0$ is a smooth curve $q:(0,\epsilon_0)\to D\cap\Lambda$ with $\lim_{t\to+0}q(t)=(0,\dots,0)$, where 
\begin{equation*}
\Lambda=\{(z_1,z')\in \C\times\C^n : \alpha|z'|<\Re z_1\}, 
\end{equation*}
for some $\alpha>0$. For a vector $\xi\in\C^{n+1}$, let $\xi_{N,t}$ and $\xi_{T,t}^i$ denote the complex normal and $i$-th complex tangential components of $\xi$ at $\pi(q(t))$, respectively. Here, $\pi:\C^{n+1}\to bD$ is the normal projection onto $bD$ defined for points sufficiently close to $bD$, so that $\xi=\xi_{N,t}+\sum_{i=1}^k\xi_{T,t}^i$, $\xi_{T,t}^i\in T_{\pi(q(t))}^{\C}(bD)$ for $i=1,\dots,k$, and $\xi_{N,t}\bot T_{\pi(q(t))}^{\C}(bD)$. Let $d(t)$ be the Euclidean distance from $q(t)$ to $bD$ along the normal direction, and for each of $i=1,\dots,k$, let $d_i^\ast(t)$ be the corresponding Euclidean distance along the $i$-th tangential direction.

The following is our main theorem:

\begin{theorem}\label{MT}
Let $D\subset\C^{n+1}$ be a smoothly bounded domain which is generally decoupled  near $0\in bD$, and let $I$ and $J$ be index sets satisfying $I\sqcup J=\{1,\dots,k\}$ such that $f_i$ in (\ref{eq12}) is finite type at $0$ for $i\in I$ and exponentially flat at $0$ for $i\in J$. If $D$ is pseudoconvex and $q$ is a nontangential cone stream approaching $0$, then, for $\xi\in\C^{n+1}\setminus\{0\}$,
\begin{equation}
\ds\lim_{t\to+0}\cal{K}_D(q(t))d(t)^{2+\sum_{i\in I}2n_i/m_i}\ds\prod_{i\in J}d_i^\ast(t)^{2n_i}=\ds\frac1{4\mathrm{vol}(\mathscr{E})\prod_{i\in J}\mathrm{vol}(B_{n_i}(0,1))},
\end{equation}
\begin{equation}
\ds\lim_{t\to+0}\frac{B_D(q(t);\xi)}{\left\{\mathrm{vol}(\mathscr{E})\left(a_1\ds\frac{|\xi_{N,t}|^2}{4d(t)^2}+\ds\sum_{i\in I}a^i\frac{|\xi_{T,t}^i|^2}{d(t)^{2/m_i}}\right)+\ds\sum_{i\in J}(n_i+1)\frac{|\xi_{T,t}^i|^2}{d_i^\ast(t)^2}\right\}^{1/2}}=1,\text{ and}
\end{equation}
\begin{equation}
\begin{split}
\ds\lim_{t\to+0}H_D(q(t);\xi)\Bigg[\mathrm{vol}(\mathscr{E})\Bigg\{\frac{b_{11}}4&\frac{|\xi_{N,t}|^4}{d(t)^4}+\sum_{i\in I}\left(b_1^i\frac{|\xi_{N,t}|^2|\xi_{T,t}^i|^2}{4d(t)^{2+2/m_i}}+4b^{ii}\frac{|\xi_{T,t}^i|^4}{d(t)^{4/m_i}}\right)\\
&+\sum_{\substack{i,j\in I\\ i<j}}b^{ij}\frac{|\xi_{T,t}^i|^2|\xi_{T,t}^j|^2}{d(t)^{2/m_i+2/m_j}}\Bigg\}+2\sum_{i\in J}(n_i+1)\frac{|\xi_{T,t}^i|^4}{d_i^\ast(t)^4}\Bigg]^{-1}\\
&\hspace{-45mm}\times\left\{\mathrm{vol}(\mathscr{E})\left(a_1\frac{|\xi_{N,t}|^2}{4d(t)^2}+\sum_{i\in I}a^i\frac{|\xi_{T,t}^i|^2}{d(t)^{2/m_i}}\right)+\sum_{i\in J}(n_i+1)\frac{|\xi_{T,t}^i|^2}{d_i^\ast(t)^2}\right\}^2=-2.
\end{split}
\end{equation}
Here, $\cal{K}_D, B_D$ and $H_D$ are the Bergman kernel of $D$ on the diagonal, the Bergman metric of $D$, and the holomorphic sectional curvature of $B_D$, respectively. The set $\mathscr{E}$ is an ellipsoid defined by
\begin{equation}\label{eq:ellipsoid}
\mathscr{E}=\left\{(z_1,z')\in\C\times\ds\prod_{i\in I}\C^{n_i} : |z_1|^2+\ds\sum_{i\in I}c_i|z^i|^{m_i}<1\right\},
\end{equation}
where $m_i\ge2$ and $c_i>0$ denote the constants given by Definition \ref{D21} (ii) for $i\in I$.  The numbers $a_1, a^i, b_{11}, b_1^i$ and $b^{ij}$ are constants determined by the reciprocals of the $L^2$-norm of certain monomials on $\mathscr{E}$ (see Lemma \ref{L26}). 
\end{theorem}

Since setting $I=\emptyset$ and $k=1$ in the above theorem yields the result of Jaiswal in \cite{Jais254}, our theorem can be seen as an extension of \cite{Jais254}. In addition, we would like to emphasize that our theorem is not a straightforward extension of previous results. This is because, in the scaling argument in Section 3, different scalings must be considered for the coexisting finite-type and exponentially flat directions. More precisely, we show that, by scaling the finite-type directions by fractional powers of the normal distance and the exponentially flat directions by the tangential distance, the resulting domain converges near the origin to the product of the ellipsoid and the unit balls. One can see from the formula in Theorem \ref{MT} that these different scalings are directly related to the boundary limits of the Bergman kernel and its associated invariants.

This paper is organized as follows. In Section 2, we present some preliminaries, including a certain class of smooth functions, and the Bergman invariants and the associated minimum integrals. In Section 3, we compute the asymptotic boundary limits of these integrals. These results are then used in Section 4 to prove Theorem \ref{MT}.

\section{Preliminaries}
In this section, we set up some of the notation and concepts for Theorem \ref{MT}. 

\subsection{Finite type and exponentially flat}
Let us introduce a class of smooth functions for the definition of generally decoupled domains. 

\begin{definition}\label{D21}
A smooth function $f:\R\to\R$ is said to be {\it finite type at the origin} if it satisfies the following conditions:
\begin{enumerate}
   \item[(i)] $f$ is strictly increasing on $[0,\infty)$, and
   \item[(ii)] there exists $m\in\N$ and a positive constant $c$ such that
   \begin{equation*}
   \ds\lim_{x\to+0}\ds\frac{f(x)}{x^m}=c.
   \end{equation*}
\end{enumerate}
From the definition, one can easily see that if $f$ is finite type at the origin, then $f(0)=0$ and $f(x)>0$ for $x>0$.
\end{definition}

\begin{example}
For $m\in\N$, $f(x)=x^m$ is finite type at the origin.
\end{example}

\begin{definition}\label{D23}
A smooth function $f:\R\to\R$ is said to be {\it exponentially flat at the origin} if it satisfies the following conditions:
\begin{enumerate}
   \item[(i)] $f(x)=0$ for $x\le0$,
   \item[(ii)] there exists a constant $\epsilon_0>0$ such that $f''(x)>0$ for $x\in(0,\epsilon_0)$, and
   \item[(iii)] for the function
   \begin{center}
   $g(x):=\left\{
   \begin{array}{cl}
    -\ds\frac1{\log{\left(f(x)\right)}} & \text{if } x\in(0,\epsilon_0),\text{ and}\\
    0 & \text{if } x=0,
   \end{array}
   \right.$
   \end{center}
   there exist positive constants $\mu$ and $C$ such that
   \begin{equation*}
   \ds\lim_{x\to+0}\ds\frac{g(x)}{x^\mu}=C.
   \end{equation*}
\end{enumerate}
\end{definition}

\begin{example}
For $\mu\in\N$, $f(x)=e^{-1/x^\mu}$ is exponentially flat at the origin.
\end{example}

\begin{remark}\label{R25}
If $f$ is exponentially flat at the origin, it is easily seen that $f$ is flat at the origin, i.e. 
\begin{equation*}
f^{(k)}(0)=0,\text{ for all }k\in\N,\text{ and }
\end{equation*}
\begin{equation*}
\ds\lim_{x\to+0}\ds\frac{f^{(k)}(x)}{x^l}=0,\text{ for all }k,l\in\N.
\end{equation*}
\end{remark}

Recall that
\begin{equation*}
I=\left\{i : f_i \text{ is finite type at the origin}\right\},
\end{equation*}
\begin{equation*}
J=\left\{i : f_i \text{ is exponentially flat at the origin}\right\},
\end{equation*}
and $I\sqcup J=\{1,\dots,k\}$. Let $m_i\in\N$ and $c_i>0$ be the constants given by Definition \ref{D21} (ii) if $i\in I$, and let $\mu_i$ and $C_i$ denote the positive constants given by Definition \ref{D23} (iii) if $i\in J$.
\begin{remark}
$m_i$ and $\mu_i$, obtained from $f_i$ in Definition \ref{D13}, are even due to the smoothness of $D$.
\end{remark}

\subsection{The Bergman invariants}
Let $D$ be a bounded domain in $\C^{n+1}$ and let $A^2(D)$ denote the Hilbert space consisting of all $L^2$-holomorphic functions on $D$, referred to as the Bergman space of $D$. Let $\|\cdot\|_D$ denote the norm in $L^2(D)$.

For $z\in D$, the Bergman kernel of $D$ on the diagonal $\cal{K}_D(z)$ is given by $\cal{K}_D(z)=\sum_\alpha|\phi_\alpha(z)|^2$, where $(\phi_\alpha)_\alpha$ is a complete orthonormal system of $A^2(D)$. The Bergman metric of $D$ is defined by
\begin{equation*}
B_D(z;\xi)=\left(\ds\sum_{j,k=1}^{n+1}g_{j\bar{k}}(z)\xi_j\bar\xi_k\right)^{1/2},\text{ for }z\in D\text{ and }\xi\in\C^{n+1},
\end{equation*}
where
\begin{equation*}
g_{j\bar{k}}(z)=\ds\frac{\partial^2\log{\cal{K}_D}}{\partial z_j\partial\bar{z}_k}(z).
\end{equation*}

The holomorphic sectional curvature of the Bergman metric is defined by 
\begin{equation*}
H_D(z;\xi)=\ds\frac{\sum_{h,j,k,l}R_{\bar{h}jk\bar{l}}(z)\bar\xi_h\xi_j\xi_k\bar\xi_l}{B_D(z;\xi)^4},\text{ for }z\in D\text{ and }\xi\in\C^{n+1}\setminus\{0\},
\end{equation*}
where
\begin{equation*}
R_{\bar{h}jk\bar{l}}(z)=-\ds\frac{\partial^2g_{j\bar{h}}}{\partial z_k\partial\bar{z}_l}(z)+\ds\sum_{\mu,\nu}g^{\nu\bar\mu}(z)\ds\frac{\partial g_{j\bar\mu}}{\partial z_k}(z)\ds\frac{\partial g_{\nu\bar{h}}}{\partial\bar{z}_l}(z),
\end{equation*}
and $(g^{\nu\bar\mu}(z))$ represents the inverse matrix of $(g_{j\bar{k}}(z))$.

\subsection{Minimum integrals}
We now define certain minimum integrals associated with $\cal{K}_D, B_D,$ and $H_D$, which are useful for computing these quantities. For a point $z\in D$ and nonzero vector $\xi\in\C^{n+1}\setminus\{0\}$, define
\begin{equation*}
I_0^D(z)=\inf{\left\{\|f\|_D^2 : f\in A^2(D), f(z)=1\right\}},
\end{equation*}
\begin{equation*}
I_1^D(z;\xi)=\inf{\left\{\|f\|_D^2 : f\in A^2(D), f(z)=0, \ds\sum_{j=1}^{n+1}\xi_j\ds\frac{\partial f}{\partial z_j}(z)=1\right\}},\text{ and }
\end{equation*}
\begin{equation*}
\begin{split}
I_2^D(z;\xi)=\inf\Bigg\{\|f\|_D^2 : f\in A^2(D), f(z)=\ds\frac{\partial f}{\partial z_1}(z)=\dots=&\ds\frac{\partial f}{\partial z_{n+1}}(z)=0,\\
&\ds\sum_{j,k=1}^{n+1}\xi_j\xi_k\ds\frac{\partial^2 f}{\partial z_j\partial z_k}(z)=1\Bigg\}.
\end{split}
\end{equation*}
From these definitions, the following monotonicity and transformation rules hold: if $D'\subset D$ is a subdomain, $z\in D'$ and $\xi\in\C^{n+1}\setminus\{0\}$, then
\begin{equation}\label{eq21}
I_0^{D'}(z)\le I_0^D(z),\text{ and }I_j^{D'}(z;\xi)\le I_j^D(z;\xi),\text{ for }j=1,2.
\end{equation}
In addition, if $F:D_1\to D_2$ is a biholomorphic map, then 
\begin{equation}\label{eq22}
I_0^{D_1}(z)\left|\mathrm{det}J_\C F(z)\right|^2=I_0^{D_2}(F(z)),\text{ and }
\end{equation}
\begin{equation}\label{eq23}
I_j^{D_1}(z;\xi)\left|\mathrm{det}J_\C F(z)\right|^2=I_j^{D_2}(F(z);J_\C F(z)\xi),\text{ for }j=1,2,
\end{equation}
where $J_\C F(z)$ denotes the complex Jacobian matrix of $F$ at $z$. It is also known that these quantities are smooth with respect to the point $z$ and the nonzero vector $\xi$, and 
\begin{equation}\label{eq24}
\cal{K}_D(z)=\ds\frac1{I_0^D(z)},\ B_D(z;\xi)^2=\ds\frac{I_0^D(z)}{I_1^D(z;\xi)},\text{ and }H_D(z;\xi)=2-\ds\frac{\left\{I_1^D(z;\xi)\right\}^2}{I_0^D(z)I_2^D(z;\xi)}.
\end{equation}
For these properties and further details, we refer the reader to, for example, \cite{JP13}. 

We write $B_d(0,r)$ for the open ball of radius $r>0$ centered at the origin in $\C^d$. The following lemma provides explicit formulas for the Bergman kernel, the Bergman metric, and its holomorphic sectional curvature on the product domain of the ellipsoid in (\ref{eq:ellipsoid}) and open balls, which are needed to prove Theorem \ref{MT}.
\begin{lemma}\label{L26}
Let $\mathscr{E}$ be the ellipsoid as in (\ref{eq:ellipsoid}). The Bergman kernel on the diagonal, the Bergman metric and its holomorphic sectional curvature at the origin in $\mathscr{E}\times\prod_{i\in J}B_{n_i}(0,1)$ are as follows.
\begin{equation}\label{eq26}
\cal{K}_{\mathscr{E}\times\prod_{i\in J}B_{n_i}(0,1)}(0)=\ds\frac1{\mathrm{vol}(\mathscr{E})\prod_{i\in J}\mathrm{vol}(B_{n_i}(0,1))},
\end{equation}
\begin{equation}\label{eq27}
B_{\mathscr{E}\times\prod_{i\in J}B_{n_i}(0,1)}(0;\xi)=\left\{\mathrm{vol}(\mathscr{E})\!\!\left(a_1|\xi_1|^2+\ds\sum_{i\in I}a^i|\xi^i|^2\right)+\ds\sum_{i\in J}(n_i+1)|\xi^i|^2\right\}^{1/2}\!\!\!\!\!\!\!,\text{ and}
\end{equation}
\begin{equation}\label{eq28}
\begin{split}
H_{\mathscr{E}\times\prod_{i\in J}B_{n_i}(0,1)}(0;\xi)=&-B_{\mathscr{E}\times\prod_{i\in J}B_{n_i}(0,1)}(0;\xi)^{-4}\\
&\times\Bigg[\mathrm{vol}(\mathscr{E})\Bigg\{4b_{11}|\xi_1|^4+\ds\sum_{i\in I}\left(b_1^i|\xi_1|^2|\xi^i|^2+4b^{ii}|\xi^i|^4\right)\\
&\hspace{22mm}+\ds\sum_{\substack{i,j\in I\\ i<j}}b^{ij}|\xi^i|^2|\xi^j|^2\Bigg\}+2\ds\sum_{i\in J}(n_i+1)|\xi^i|^4\Bigg],
\end{split}
\end{equation}
where $a_1=\|z_1\|_{\mathscr{E}}^{-2}, a^i=\|z_1^i\|_{\mathscr{E}}^{-2}, b_{11}=\|z_1z_1\|_{\mathscr{E}}^{-2}, b_1^i=\|z_1z_1^i\|_{\mathscr{E}}^{-2}$ and $b^{ij}=\|z_1^iz_1^j\|_{\mathscr{E}}^{-2}$ for $i,j\in I$. Note that, due to the circular symmetries of $\mathscr{E}$, these values are independent of the choice of the components $z_k^i$ $(k=1,\dots,n_i)$; thus, the first components are chosen as representatives.
\end{lemma}

\begin{remark}\label{R27}
$\mathscr{E}$ in (\ref{eq:ellipsoid}) is a bounded domain biholomorphically equivalent, via the map $\Psi$ given by
\begin{equation*}
\Psi_1(z_1,z')=\ds\frac{1-z_1}{1+z_1},\ \Psi^i(z_1,z')=\ds\frac{z^i}{(1+z_1)^{2/m_i}},\text{ for }i\in I,
\end{equation*}
to the unbounded domain
\begin{equation*}
\mathscr{U}=\left\{(z_1,z')\in\C\times\ds\prod_{i\in I}\C^{n_i} : \ds\sum_{i\in I}c_i|z^i|^{m_i}<\Re z_1\right\},
\end{equation*}
called the polynomial half space. The Bergman kernels of these domains have been investigated by, for instance, Haslinger \cite{Has94, Has98}, Chakrabarti and Upadrashta \cite{CU19}, and Kamimoto \cite{Kami04, Kami24}.
\end{remark}

\begin{proof}[Proof of Lemma \ref{L26}]
 Let $z\in\mathscr{E}$ and $\xi\in\C^{n+1}\setminus\{0\}$. Since the domain $\mathscr{E}$ has circular symmetries, an orthonormal system for $A^2(\mathscr{E})$ is given by the monomials, which implies 
\begin{equation}\label{eq29}
\cal{K}_{\mathscr{E}}(z)=\ds\sum_{\alpha}\ds\frac{|z^\alpha|^2}{\|z^\alpha\|_{\mathscr{E}}^2},
\end{equation}
where
\begin{equation*}
|z^\alpha|^2=|z_1|^{2\alpha_1}\prod_{i\in I}|z^i|^{2\alpha^i},\text{ for }\alpha\in\Z_+\times\prod_{i\in I}\Z_+^{n_i},\text{ and }
\end{equation*}
\begin{equation*}
|z^i|^{2\alpha^i}=|z_1^i|^{2\alpha_1^i}\cdots|z_{n_i}^i|^{2\alpha_{n_i}^i}.
\end{equation*}
By the product formula for the Bergman kernel on a product domain, we get (\ref{eq26}).

Now recall that the Bergman kernel of a unit ball in $\C^{n_i}$ on the diagonal is given by
\begin{equation}\label{eq210}
\cal{K}_{B_{n_i}(0,1)}(z^i)=\ds\frac{n_i!}{\pi^{n_i}}\ds\frac1{(1-|z^i|^2)^{n_i+1}}.
\end{equation}

Considering (\ref{eq29}) and (\ref{eq210}), we have 
\begin{equation*}
B_{\mathscr{E}}(0;\xi)^2=\mathrm{vol}(\mathscr{E})\left(a_1|\xi_1|^2+\ds\sum_{i\in I}a^i|\xi^i|^2\right),\text{ and}
\end{equation*}
\begin{equation*}
B_{B_{n_i}(0,1)}(0;\xi)^2=(n_i+1)|\xi^i|^2,\text{ for }i\in J,
\end{equation*}
where $a_1=\|z_1\|_{\mathscr{E}}^{-2}, a^i=\|z_1^i\|_{\mathscr{E}}^{-2}$, which implies (\ref{eq27}).

Similarly, using (\ref{eq29}) and (\ref{eq210}) again, we obtain (\ref{eq28}).
\end{proof}

\section{Boundary behavior of minimum integrals}
Let $D$ be the domain as in Theorem \ref{MT}, and let $q$ be a nontangential cone stream approaching $0\in bD$. In this section, we will study the boundary behavior of the minimum integrals of $D$ along $q$, which is used in the proof of Theorem \ref{MT} directly. Since $D$ is generally decoupled near $0$, we get an open neighborhood $U$ of $0$ and a local defining function $\rho$ of $D$ near $0$ as in Definition \ref{D13} such that
\begin{equation*}
D\cap U=\left\{(z_1,z')\in U : \rho(z_1,z')<0\right\}.
\end{equation*}

To prove Theorem \ref{MT}, we first change the coordinates so that the stream $q$ becomes a straight line approaching $0$ along the $\Re z_1$-axis, that is, in the normal direction. In Subsection 3.2, we suitably cut out a window-shaped domain of $D$ near $0$ in the new coordinates and apply a scaling map and a biholomorphic map to this domain, which implies that it converges to $\mathscr{E}\times\prod_{i\in J}B_{n_i}(0,1)$. Note that the scaling map we construct in this subsection applies different scalings along the finite-type and exponentially flat directions. At the end of Subsection 3.2, we show that the minimum integrals of $D$ and of the window-shaped domain in the new coordinates exhibit the same behavior as $t\to+0$. For this localization, we will only state the conclusion because the proof is the same as that of Lemma 3.5 in \cite{Jais254}. Finally, in Subsection 3.3, we compute the precise boundary limits of the minimum integrals of $D$.

\subsection{Centering}
For each $t>0$, we first apply translations and unitary maps to the domain $D$ so that
\begin{equation*}
q(t)\mapsto(\text{dist}(q(t),bD),0,\dots,0). 
\end{equation*}
For $t>0$, define the translation map $T_t^1:\C^{n+1}\to\C^{n+1}$ by
\begin{equation*}
T_t^1(w):=\left(w_1-i\Im q_1(t),w'\right).
\end{equation*}
So, $T_t^1(q(t))=\left(\Re q_1(t),q^1(t),\dots,q^k(t)\right)$. 

We next choose the unitary map so that
\begin{equation*}
\left(\Re q_1(t),q^1(t),\dots,q^k(t)\right)\mapsto\left(\Re q_1(t),\left(\|q^1(t)\|_{\C^{n_1}},0,\dots,0\right),\dots,\left(\|q^k(t)\|_{\C^{n_k}},0,\dots,0\right)\right).
\end{equation*}
For the purpose, set $r_i(t):=\|q^i(t)\|_{\C^{n_i}}$ for each $i=1,\dots,k$. In each $\C^{n_i}$, we choose a unitary transformation $R_{t,i}:\C^{n_i}\to\C^{n_i}$ that rotates $q^i(t)$ onto the real first coordinate axis so that
\begin{equation*}
R_{t,i}(x)=\left\{
   \begin{array}{cl}
   x & \text{ if } r_i(t)=0,\text{ and}\\
   e_1^i & \text{ if } x=\ds\frac{q^i(t)}{r_i(t)}\ \text{ and }\ r_i(t)\ne0,
   \end{array}
   \right.
\end{equation*}
where $e_1^i=(1,0,\dots,0)\in\C^{n_i}$, from which it follows that 
\begin{equation*}
R_{t,i}(q^i(t))=\left\{
   \begin{array}{cl}
   0 & \text{ if } q^i(t)=0,\text{ and}\\
   (r_i(t),0,\dots,0) & \text{ if } q^i(t)\ne0.
   \end{array}
   \right.
\end{equation*}
Now define the unitary map $R_t^1:\C^{n+1}\to\C^{n+1}$ by
\begin{equation*}
R_t^1=\mathrm{id}_\C\oplus R_{t,1}\oplus\dots\oplus R_{t,k}.
\end{equation*}

\begin{remark}\label{R31}
A domain $\Omega:=\left\{z\in\C^{n+1} : \rho(z)<0\right\}$ is invariant under both $T_t^1$ and $R_t^1$ for $t>0$ and
\begin{equation*}
\Omega\cap U=D\cap U.
\end{equation*}
\end{remark}

We next set
\begin{equation*}
\tilde{q}(t):=R_t^1\circ T_t^1(q(t))=\left(\Re q_1(t),\left(r_1(t),0,\dots,0\right),\dots,\left(r_k(t),0,\dots,0\right)\right).
\end{equation*}
Let $\delta_0>0$ such that $(-4\delta_0,4\delta_0)^2\times\prod_{i=1}^kB_{n_i}(0,4\delta_0)\subset U$. Since $q(t)\to0$ as $t\to+0$, $q(t)\in(-\delta_0,\delta_0)^2\times\prod_{i=1}^kB_{n_i}(0,\delta_0)\subset U$ for sufficiently small $t>0$. Hence, for any sufficiently small $t>0$, we get
\begin{equation*}
\tilde{q}(t)\in\Omega\cap R_t^1\circ T_t^1\left((-\delta_0,\delta_0)^2\times\ds\prod_{i=1}^kB_{n_i}(0,\delta_0)\right)\cap\Lambda\subset\Omega\cap U\cap\Lambda,
\end{equation*}
where $\Lambda=\{(z_1,z') : \alpha|z'|<\Re z_1\}$ for some $\alpha>0$. For small $t>0$, there exists a unique point
\begin{equation*}
p(t)=\left(p_1(t),\left(p_1^1(t),0,\dots,0\right),\dots,\left(p_1^k(t),0,\dots,0\right)\right)\in b\Omega
\end{equation*}
such that $d(t):=\text{dist}(q(t),b\Omega)=\text{dist}(\tilde{q}(t),b\Omega)=|\tilde{q}(t)-p(t)|, p_1(t)\in\R, p_1^i(t)\ge0$ for $i=1,\dots,k$ and $p_1(t)=\sum_{i=1}^kf_i(p_1^i(t))\ge0$.

Now we choose another translation which maps $p(t)$ to the origin, and a unitary transformation which sends the outer normal vector of $\Omega$ at $p(t)$ to the direction $-e_1$. For $t>0$, define the translation map $T_t^2:\C^{n+1}\to\C^{n+1}$ by
\begin{equation*}
T_t^2(w):=w-p(t).
\end{equation*}
Noting that $\nabla\rho(p(t))=\left(-1,(f_1'(p_1^1(t)),0,\dots,0),\dots,(f_k'(p_1^k(t)),0,\dots,0)\right)$, the desired unitary transformation is given as follows. Define a unitary matrix $R_t'$, acting on $(z_1,z_1^1,\dots,z_1^k)^\top\in\C^{k+1}$, by
\begin{equation*}
R_t':=\ds\begin{pmatrix}
   \frac1{A(t)} & -\frac{b_1(t)}{A(t)} & -\frac{b_2(t)}{A(t)} & \cdots & -\frac{b_k(t)}{A(t)} \\
   \frac{b_1(t)}{A(t)} & 1-D_{11}(t) & -D_{12}(t) & \cdots & -D_{1k}(t) \\
   \frac{b_2(t)}{A(t)} & -D_{21}(t) & 1-D_{22}(t) &\cdots & -D_{2k}(t) \\
   \vdots & \vdots & \vdots & \ddots & \vdots \\
   \frac{b_k(t)}{A(t)} & -D_{k1}(t) & -D_{k2}(t) & \cdots & 1-D_{kk}(t)
   \end{pmatrix},
\end{equation*}
where $A(t):=\|\nabla\rho(p(t))\|, b_i(t)=f_i'(p_1^i(t))$ and $D_{ij}(t)=(b_i(t)b_j(t))/(A(t)(A(t)+1))$. We then define $R_t^2:\C^{n+1}\to\C^{n+1}$ to be the unitary transformation given by the direct sum of $R_t'$ and the identical maps on the remaining coordinates $(z_2^i,\dots,z_{n_i}^i), i=1,\dots,k$. Indeed, as $A(t)^2=1+\sum_{i=1}^kb_i(t)^2$ and $b_i(t)\in\R$, a direct computation shows that the above matrix is unitary. Furthermore, as $p(t)=\tilde{q}(t)+(\nabla\rho(p(t))/A(t))d(t)$, one can readily check that $R_t^2\circ T_t^2(\tilde{q}(t))=(d(t),0,\dots,0)$.

For $t>0$, we set $\gamma_t=R_t^2\circ T_t^2$. Letting $Z=\gamma_t^{-1}(z)$ for $z=(z_1,z^1,\dots,z^k)\in\C^{n+1}$, we have
\begin{equation*}
Z_1=\ds\frac{z_1}{A(t)}+p_1(t)+\ds\sum_{i=1}^k\ds\frac{b_i(t)}{A(t)}z_1^i, 
\end{equation*}
\begin{equation*}
Z^i=\left(-\ds\frac{b_i(t)}{A(t)}z_1+z_1^i+p_1^i(t)-b_i(t)\ds\sum_{j=1}^k\ds\frac{b_j(t)}{A(t)(A(t)+1)}z_1^j,z_2^i,\dots,z_{n_i}^i\right)
\end{equation*}
for $i=1,\dots,k$, and 
\begin{equation*}
\rho\circ\gamma_t^{-1}(z)=\ds\sum_{i=1}^kf_i(|Z^i|)-\Re Z_1.
\end{equation*}

Set
\begin{equation*}
W:=\left(-\ds\frac{\delta_0}{5(k+2)},\ds\frac{\delta_0}{5(k+2)}\right)^2\times\ds\prod_{i=1}^kB_{n_i}\left(0,\ds\frac{\delta_0}{5(k+2)}\right).
\end{equation*}
Then, for all sufficiently small $t>0$, $W\subset\gamma_t\circ R_t^1\circ T_t^1(U)$, which implies
\begin{equation*}
\gamma_t\left(R_t^1\circ T_t^1(D)\right)\cap W=\gamma_t(\Omega)\cap W
\end{equation*}
for all sufficiently small $t>0$, by combining with the relation of sets in Remark \ref{R31}. In other words, after applying $T_t^1, R_t^1$, and $\gamma_t$ to the domain $D$ and $\gamma_t$ to the domain $\Omega$, the resulting domains coincide in a small neighborhood of $0$.

\subsection{Cutting, scaling and localization}
In this subsection, we cut out the domain $\gamma_t(\Omega)$ near $0$ and apply a scaling map together with a biholomorphic map to the resulting domain so that it converges to $\mathscr{E}\times\prod_{i\in J}B_{n_i}(0,1)$. 

We begin by defining a window-shaped domain near the origin in $\gamma_t(\Omega)$, the tangential distances in each coordinate direction, and the scaling map applying to this domain. For $\epsilon,t>0$, define
\begin{equation*}
D_t^\epsilon:=\left\{(z_1,z^1,\dots,z^k)\in W : \rho\circ\gamma_t^{-1}(z_1,z^1,\dots,z^k)<0, \Re z_1<d(t)^{\frac1{(1+\epsilon)^2}}\right\},
\end{equation*}
and define
\begin{equation*}
d_i^\ast(t):=\min{\left\{s>0 : (d(t),0,\dots,0)+se_1^i\in b\gamma_t(\Omega)\right\}}\text{ and }
\end{equation*}
\begin{equation*}
\begin{split}
d_i^\epsilon(t):=\sup \Bigg\{|w^i| : (w_1,0,\dots,w^i,\dots,&0)\in W,\\
&\ds\sum_{j=1}^kf_j(|\zeta_{(i)}^j|)<\ds\frac{d(t)^{\frac1{(1+\epsilon)^2}}}{A(t)}+p_1(t)+\ds\frac{b_j(t)}{A(t)}\Re w_1^j\Bigg\},
\end{split}
\end{equation*}
for each $i\in J$. Here, $\zeta_{(i)}=(\zeta_1,\zeta_{(i)}^1,\dots,\zeta_{(i)}^k):=\gamma_t^{-1}(w_1,0,\dots,w^i,\dots,0)$ for $i=1,\dots,k$.
Define a scaling map $\Sigma_t=(\Sigma_{t,1},\Sigma_t^1,\dots,\Sigma_t^k):\C^{n+1}\to\C^{n+1}$ given by
\begin{equation*}
\Sigma_{t,1}(z_1,z')=\ds\frac{z_1}{d(t)},\ \Sigma_t^i(z_1,z')=\left\{
   \begin{array}{cl}
    \ds\frac{z^i}{d(t)^{1/m_i}} & ,\text{ if }i\in I,\text{ and} \\
    \ds\frac{z^i}{d_i^\ast(t)} & ,\text{ if }i\in J
   \end{array}
   \right.
\end{equation*}
for $t>0$. Notice that $\Sigma_t$ applies the algebraic scaling $d(t)^{1/m_i}$ to the finite-type coordinates ($i\in I$), while simultaneously applying the logarithmic tangential scaling $d_i^\ast(t)$ to the exponentially flat coordinates ($i\in J$). This multi-scale transformation yields the product of the ellipsoid and unit balls as the limit domain (see Lemma \ref{L33} below).

We first establish the asymptotic relations, as $t\to+0$, among the quantities $p(t)$, $d(t)$, $d_i^\ast(t)$ and $d_i^\epsilon(t)$.

\begin{lemma}\label{L32}
Let $D\subset\C^{n+1}$ and $U\subset\C^{n+1}$ be as defined at the beginning of Section 3. Then, for $i=1,\dots,k$ and $j\in J$, 
\begin{equation*}
\ds\lim_{t\to+0}\ds\frac{p_1^i(t)}{d(t)}\le \ds\frac k\alpha,\ \ \ \ds\lim_{t\to+0}\ds\frac{p_1^i(t)}{d_j^\ast(t)}=0,\text{ and }\ds\liminf_{t\to+0}\ds\frac{d_j^\ast(t)}{d_j^\epsilon(t)}\ge\ds\frac1{(1+\epsilon)^{2/\mu_j}}.
\end{equation*}
Here, recall that $\alpha>0$ is the constant determining the aperture of the cone $\Lambda$, which appears in the definition of the nontangential cone stream (see the notation introduced just before Theorem \ref{MT}).
\end{lemma}

\begin{proof}
Without loss of generality, we may assume that each $p_1^i(t)$ is nonzero for all small $t>0$. Note that, from Definition \ref{D21} and Remark \ref{R25},
\begin{equation*}
\ds\lim_{t\to+0}\ds\frac{f_i(p_1^i(t))}{p_1^i(t)}=0,\ \text{ for }i=1,\dots,k,
\end{equation*}
\begin{equation*}
\ds\lim_{t\to+0}\ds\frac{b_i(t)}{p_1^i(t)}\le c_im_i,\ \text{ for }i\in I,\text{ and }\ds\lim_{t\to+0}\ds\frac{b_i(t)}{p_1^i(t)}=0,\ \text{ for }i\in J.
\end{equation*}
As $\tilde{q}(t)\in\Lambda$ and $p(t)=\tilde{q}(t)+\frac{\nabla\rho(p(t))}{A(t)}d(t)$, we have
\begin{equation}\label{eq31}
\alpha\left(\ds\sum_{i=1}^kr_i(t)^2\right)^{1/2}<\Re q_1(t),
\end{equation}
\begin{equation}\label{eq32}
p_1(t)=\Re q_1(t)-\ds\frac{d(t)}{A(t)},\text{ and }
\end{equation}
\begin{equation}\label{eq33}
p_1^i(t)=r_i(t)+\ds\frac{b_i(t)}{A(t)}d(t),\text{ for }i=1,\dots,k.
\end{equation}
Dividing both side of (\ref{eq33}) by $p_1^i(t)$, we have
\begin{equation*}
1=\ds\frac{r_i(t)}{p_1^i(t)}+\ds\frac{b_i(t)}{p_1^i(t)}\ds\frac{d(t)}{A(t)},
\end{equation*}
which implies
\begin{equation*}
\ds\lim_{t\to+0}\ds\frac{r_i(t)}{p_1^i(t)}=1,\ \ \text{for }i=1,\dots,k.
\end{equation*}
Since, from (\ref{eq31}), (\ref{eq32}) and $p_1(t)=\sum_{i=1}^kf_i(p_1^i(t))$,
\begin{align*}
\ds\sum_{i=1}^kf_i(p_1^i(t))&>\alpha\left(\ds\sum_{i=1}^kr_i(t)^2\right)^{1/2}-\ds\frac{d(t)}{A(t)}\\
&\ge\ds\frac\alpha k\ds\sum_{i=1}^kr_i(t)-d(t),
\end{align*}
we get
\begin{equation}\label{eq34}
d(t)>\ds\sum_{i=1}^k\left\{\ds\frac\alpha kr_i(t)-f_i(p_1^i(t))\right\}.
\end{equation}
Considering 
\begin{equation*}
\ds\frac\alpha kr_i(t)-f_i(p_1^i(t))=p_1^i(t)\left[\ds\frac\alpha k\frac{r_i(t)}{p_1^i(t)}-\frac{f_i(p_1^i(t))}{p_1^i(t)}\right],
\end{equation*}
it follows that $\frac\alpha kr_i(t)-f_i(p_1^i(t))\ge0$ for sufficiently small $t>0$, and hence, together with (\ref{eq34}), we have
\begin{equation*}
\ds\frac{p_1^i(t)}{d(t)}<\ds\frac1{\frac\alpha k\frac{r_i(t)}{p_1^i(t)}-\frac{f_i(p_1^i(t))}{p_1^i(t)}},
\end{equation*}
which implies that $\lim_{t\to+0}\frac{p_1^i(t)}{d(t)}\le\frac k\alpha$ for $i=1,\dots,k$.

Let $i=1,\dots,k$ and $j\in J$. As $(d(t),0,\dots,d_j^\ast(t),\dots,0)\in b\gamma_t(\Omega)$, 
\begin{equation}\label{eq35}
\begin{split}
\ds\frac{d(t)}{A(t)}+\frac{b_j(t)}{A(t)}&d_j^\ast(t)+p_1(t)\\
&=f_j\left(\left|-\ds\frac{b_j(t)}{A(t)}d(t)+p_1^j(t)+\left[1-\ds\frac{b_j(t)^2}{A(t)(A(t)+1)}\right]d_j^\ast(t)\right|\right) \\
&\ \ +\ds\sum_{m\ne j}f_m\left(\left|-\ds\frac{b_m(t)}{A(t)}d(t)+p_1^m(t)-\ds\frac{b_j(t)b_m(t)}{A(t)(A(t)+1)}d_j^\ast(t)\right|\right).
\end{split}
\end{equation}
Since each $f_i$ is strictly increasing near the origin and $|-\frac{b_i(t)}{A(t)}d(t)+p_1^i(t)|\le p_1^i(t)$ for sufficiently small $t>0$, it follows, together with (\ref{eq35}), that
\begin{align*}
\ds\frac{d(t)}{A(t)}+\ds\frac{b_j(t)}{A(t)}d_j^\ast(t)&\le\ds\sum_{i=1}^k\left[f_i(d_j^\ast(t)+p_1^i(t))-f_i(p_1^i(t))\right] \\
&=\left(\ds\sum_{i=1}^kb_i(t)\right)d_j^\ast(t)+o(d_j^\ast(t))
\end{align*}
as $t\to+0$, which implies that $d(t)/d_j^\ast(t)\to0$ as $t\to+0$. Therefore, we obtain
\begin{equation*}
\ds\lim_{t\to+0}\ds\frac{p_1^i(t)}{d_j^\ast(t)}=\ds\lim_{t\to+0}\ds\frac{p_1^i(t)}{d(t)}\frac{d(t)}{d_j^\ast(t)}=0.
\end{equation*}

Finally, we prove the last inequality in the statement. By using (\ref{eq35}) again, we have $d_j^\ast(t)\sim f_j^{-1}(d(t))$ as $t\to+0$, that is, from Definition \ref{D23}, $d_j^\ast(t)\sim(-C_j\log{d(t)})^{-1/\mu_j}$ as $t\to+0$. Also, from the inequality in the definition of $d_i^\epsilon(t)$, it follows that $d_j^\epsilon(t)\lesssim(1+\epsilon)^{2/\mu_j}(-C_j\log{d(t)})^{-1/\mu_j}$ as $t\to+0$, and thus we obtain the last inequality in the lemma.
\end{proof}

Let $\Phi:\left\{z\in\C^{n+1} : \Re z_1>0\right\}\to\D\times\C^n$ denote a biholomorphic map given by
\begin{equation*}
\Phi_1(z)=\ds\frac{1-z_1}{1+z_1},\ \ \Phi^i(z)=\left\{
   \begin{array}{cl}
    \left(\ds\frac2{1+z_1}\right)^{2/m_i}z^i & \text{ if }i\in I,\text{ and}\\
    z^i & \text{ if }i\in J,
   \end{array}
   \right.
\end{equation*}
where $\D$ is the unit disc in $\C$. Note that $\Phi_1=\Psi^{-1}$ and $\Phi^i=(\Psi^i)^{-1}$ for $i\in I$, where $\Psi_1, \Psi^i$ are given in Remark \ref{R27}.

For sufficiently small $\epsilon, \delta>0$, define
\begin{equation*}
\mathscr{E}_\delta:=\left\{(z_1,z')\in\C\times\ds\prod_{i\in I}\C^{n_i} : |z_1|^2+\ds\sum_{i\in I}c_i|z^i|^{m_i}<(1-\delta)^2\right\},\text{ and}
\end{equation*}
\begin{equation*}
\mathscr{E}^\epsilon:=\left\{(z_1,z')\in\C\times\ds\prod_{i\in I}\C^{n_i} : |z_1|^2+\ds\sum_{i\in I}(c_i-\epsilon)|z^i|^{m_i}<1\right\}.
\end{equation*}
The following lemma shows that the limit domain of $\Phi\circ\Sigma_t(D_t^\epsilon)$ as $t\to+0$ is $\mathscr{E}\times\prod_{i\in J}B_{n_i}(0,1)$.

\begin{lemma}\label{L33}
Let $D\subset\C^{n+1}$ and $U\subset\C^{n+1}$ be as defined at the beginning of Section 3. Then, for every sufficiently small $\epsilon,\delta>0$, there exists $t_0(\epsilon,\delta)>0$ such that
\begin{equation*}
\mathscr{E}_\delta\times\ds\prod_{i\in J}B_{n_i}(0,1-\delta)\subset\Phi\circ\Sigma_t(D_t^\epsilon)\subset\mathscr{E}^\epsilon\times\ds\prod_{i\in J}B_{n_i}\left(0,\ds\frac{d_i^\epsilon(t)}{d_i^\ast(t)}\right)
\end{equation*}
for all $t\in(0,t_0(\epsilon,\delta))$.
\end{lemma}

\begin{proof}
Fix sufficiently small $\epsilon,\delta>0$. We first prove the first inclusion in the statement. For the purpose of this, let $z\in\mathscr{E}_\delta\times\prod_{i\in J}B_{n_i}(0,1-\delta)$. Setting $w=\Sigma_t^{-1}\circ\Phi^{-1}(z)$ yields 
\begin{equation*}
w_1=d(t)\ds\frac{1-z_1}{1+z_1},\ \ w^i=\left\{
   \begin{array}{cl}
    \left(\ds\frac{d(t)^{1/2}}{1+z_1}\right)^{2/m_i}z^i & \text{ if }i\in I,\text{ and}\\
    d_i^\ast(t)z^i & \text{ if }i\in J,
   \end{array}
   \right.
\end{equation*}
and thus it suffices to check that $w\in D_t^\epsilon$ for sufficiently small $t>0$. Noting that there exists $\zeta\in\mathscr{E}\times\prod_{i\in J}B_{n_i}(0,1)$ such that 
\begin{equation*}
z_1=(1-\delta)\zeta_1,\ \ z^i=\left\{
   \begin{array}{cl}
    (1-\delta)^{2/m_i}\zeta^i & \text{ if }i\in I,\text{ and}\\
    (1-\delta)\zeta^i & \text{ if }i\in J,
   \end{array}
   \right.
\end{equation*}
we see that
\begin{equation*}
\Re w_1=d(t)\Re\left(\ds\frac{1-(1-\delta)\zeta_1}{1+(1-\delta)\zeta_1}\right)<\ds\frac{2-\delta}\delta d(t)^{\frac\epsilon{1+\epsilon}}d(t)^{\frac1{1+\epsilon}}.
\end{equation*}
Since $d(t)^{\frac\epsilon{1+\epsilon}}\to0$ and $d(t)^{\frac1{1+\epsilon}}<d(t)^{\frac1{(1+\epsilon)^2}}$ as $t\to+0$, there exists small $t_0(\epsilon,\delta)>0$ such that
\begin{equation}\label{eq36}
\Re w_1<d(t)^{\frac1{(1+\epsilon)^2}},\ \text{ for all }t\in(0,t_0(\epsilon,\delta)).
\end{equation}
Also, as $d_i^\ast(t)\to0$ as $t\to+0$ for each $i\in J$, we can ensure that $w\in W$ for all $t\in(0,t_0(\epsilon,\delta))$ by taking $t_0(\epsilon,\delta)>0$ smaller if necessary. We next show that $\rho\circ\gamma_t^{-1}(w)<0$ for sufficiently small $t>0$. Letting $Z=\gamma_t^{-1}(w)$, $\rho\circ\gamma_t^{-1}(w)$ can be written in the form
\begin{equation*}
\rho\circ\gamma_t^{-1}(w)=\ds\sum_{i=1}^k\left[f_i(|Z^i|)-f_i(p_1^i(t))-\ds\frac{b_i(t)}{A(t)}\Re w_1^i\right]-\ds\frac{\Re w_1}{A(t)},
\end{equation*}
where 
\begin{equation*}
Z^i=\left(-\ds\frac{b_i(t)}{A(t)}w_1+p_1^i(t)+w_1^i-\ds\frac{b_i(t)}{A(t)(A(t)+1)}\ds\sum_{j=1}^kb_j(t)w_1^j,w_2^i,\dots,w_{n_i}^i\right)
\end{equation*}
for each $i=1,\dots,k$. By the definition of $\mathscr{E}_\delta$, we have
\begin{align*}
\Re w_1&>\ds\frac\delta{2-\delta}d(t)+\ds\frac{d(t)}{|1+z_1|^2}\ds\sum_{i\in I}c_i|z^i|^{m_i}\\
&=\ds\frac\delta{2-\delta}d(t)+\ds\sum_{i\in I}c_i|w^i|^{m_i},
\end{align*}
where we have used (\ref{eq36}) in the second equality. For $i\in I$, Lemma \ref{L32} together with $m_i\ge2$ yields $p_1^i(t)=o(d(t)^{1/m_i})$ and $b_i(t)=o(d(t)^{1-1/m_i})$ as $t\to+0$. Furthermore, it follows from (\ref{eq36}) and our choice of $\epsilon$ that $w_1=o(d(t)^{1/m_i})$ as $t\to+0$. Hence, we get 
\begin{equation*}
f_i(|Z^i|)-f_i(p_1^i(t))=c_i|w^i|^{m_i}+o(d(t)),\ \ \text{as }t\to+0.
\end{equation*}
Noting that $\Re w_1^i=O(d(t)^{1/m_i})$, together with $b_i(t)=o(d(t)^{1-1/m_i})$ and the above equation, we have
\begin{equation*}
\ds\sum_{i\in I}\left[f_i(|Z^i|)-f_i(p_1^i(t))-\ds\frac{b_i(t)}{A(t)}\Re w_1^i\right]=\ds\sum_{i\in I}c_i|w^i|^{m_i}+o(d(t))
\end{equation*}
as $t\to+0$. In the case where $i\in J$, one can easily see, from properties of exponentially flat functions, that
\begin{equation*}
\ds\sum_{i\in J}\left[f_i(|Z^i|)-f_i(p_1^i(t))-\ds\frac{b_i(t)}{A(t)}\Re w_1^i\right]=o(d(t)),\ \ \text{as }t\to+0.
\end{equation*}
Hence, it follows that
\begin{equation}\label{eq37}
\rho\circ\gamma_t^{-1}(w)<-\ds\frac\delta{2-\delta}\ds\frac{d(t)}{A(t)}+\frac{A(t)-1}{A(t)}\ds\sum_{i\in I}c_i|w^i|^{m_i}+o(d(t)),\ \ \text{as }t\to+0.
\end{equation}
Also, as $\frac{A(t)-1}{A(t)}=o(1)$ and $\sum_{i\in I}c_i|w^i|^{m_i}=O(d(t))$ as $t\to+0$,  (\ref{eq37}) implies that
\begin{equation*}
\rho\circ\gamma_t^{-1}(w)<-\ds\frac\delta{2-\delta}\frac{d(t)}{A(t)}+o(d(t)),
\end{equation*}
and therefore we obtain that $\rho\circ\gamma_t^{-1}(w)<0$ for all sufficiently small $t>0$.

To prove the second inclusion in the statement, let $u\in\Phi\circ\Sigma_t(D_t^\epsilon)$. Then, there exists $w\in D_t^\epsilon$ such that 
\begin{equation*}
u_1=\ds\frac{1-w_1/d(t)}{1+w_1/d(t)},\ \ u^i=\left\{
  \begin{array}{cl}
   \left(\ds\frac2{1+w_1/d(t)}\right)^{2/m_i}\ds\frac{w^i}{d(t)^{1/m_i}} & \text{ if }i\in I,\text{ and}\\
   \ds\frac{w^i}{d_i^\ast(t)} & \text{ if }i\in J.
  \end{array}
  \right.
\end{equation*}
The same argument as the above proof shows that
\begin{equation*}
\Re w_1>\ds\sum_{i\in I}(c_i-\epsilon)|w^i|^{m_i}+o(d(t)), \ \ \text{as }t\to+0.
\end{equation*}
Consequently, for sufficiently small $t>0$, we have
\begin{align*}
|u_1|^2+\ds\sum_{i\in I}(c_i-\epsilon)|u^i|^{m_i}&=\left|\ds\frac{1-w_1/d(t)}{1+w_1/d(t)}\right|^2+\ds\frac1{|1+w_1/d(t)|^2d(t)}\ds\sum_{i\in I}(c_i-\epsilon)|w^i|^{m_i}\\
&<\ds\frac1{|1+w_1/d(t)|^2}\left[|1-w_1/d(t)|^2+\ds\frac{\Re w_1}{d(t)}\right]\\
&=\ds\frac1{|1+w_1/d(t)|^2}\cdot|1+w_1/d(t)|^2\\
&=1,
\end{align*}
and thus we obtain that $(u_1,u_{\text{finite}})\in\mathscr{E}^\epsilon$ for all sufficiently small $t>0$, where $u_{\text{finite}}=(u^i)_{i\in I}$. Also, it follows from the definition of $d_i^\epsilon(t)$ that
\begin{equation*}
|u^i|=\ds\frac{|w^i|}{d_i^\ast(t)}<\ds\frac{d_i^\epsilon(t)}{d_i^\ast(t)},
\end{equation*}
which implies that $u^i\in B_{n_i}(0,d_i^\epsilon(t)/d_i^\ast(t))$ for $i\in J$. Therefore, this completes the proof of the second inclusion in the statement.
\end{proof}

After establishing the scaling lemma, we turn to the localization near $0$. Since the following localization lemma on the minimum integrals can be proved by exactly the same argument as in the proof of Lemma 3.5 in \cite{Jais254}, we only state the results.

\begin{lemma}\label{L34}
Let $D\subset\C^{n+1}$ and $U\subset\C^{n+1}$ be as defined at the beginning of Section 3 and let $D_t=\gamma_t(R_t^1\circ T_t^1(D))$. Then, for small $\epsilon>0$, 
\begin{equation*}
\ds\lim_{t\to+0}\ds\frac{I_0^{D_t}(d(t),0')}{I_0^{D_t^\epsilon}(d(t),0')}=1,\text{ and }\ds\lim_{t\to+0}\ds\frac{I_j^{D_t}((d(t),0');\zeta(t))}{I_j^{D_t^\epsilon}((d(t),0');\zeta(t))}=1
\end{equation*}
\end{lemma}
for $j=1,2$ and $\zeta(t)\in\C^{n+1}\setminus\{0\}$.

\subsection{Boundary behavior of minimum integrals}
We now compute the boundary limits of the minimum integrals of $D_t$ along $q$ by use of Lemma \ref{L33} and \ref{L34}.

\begin{lemma}\label{L35}
Let $D\subset\C^{n+1}, q$ and $U\subset\C^{n+1}$ be as defined at the beginning of Section 3, and let $D_t=\gamma_t(R_t^1\circ T_t^1(D))$. Then 
\begin{equation*}
\ds\lim_{t\to+0}\ds\frac{4d(t)^{2+\sum_{i\in I}2n_i/m_i}\prod_{i\in J}d_i^\ast(t)^{2n_i}I_0^{\mathscr{E}\times\prod_{i\in J}B_{n_i}(0,1)}(0)}{I_0^{D_t}(d(t),0')}=1,\text{ and }
\end{equation*}
\begin{equation*}
\ds\lim_{t\to+0}\ds\frac{4d(t)^{2+\sum_{i\in I}2n_i/m_i}\prod_{i\in J}d_i^\ast(t)^{2n_i}I_j^{\mathscr{E}\times\prod_{i\in J}B_{n_i}(0,1)}(0;(\Phi\circ\Sigma_t)'(d(t),0')\xi(t))}{I_j^{D_t}((d(t),0');\xi(t))}=1
\end{equation*}
for $j=1,2$ and $\xi(t)\in\C^{n+1}\setminus\{0\}$.
\end{lemma}

\begin{proof}
Let $t>0$ be sufficiently small and let $\xi(t)\in\C^{n+1}\setminus\{0\}$. Since
\begin{equation*}
\left|\mathrm{det}J_\C\Phi\circ\Sigma_t(d(t),0')\right|^2=\ds\frac1{4d(t)^{2+\sum_{i\in I}2n_i/m_i}\prod_{i\in J}d_i^\ast(t)^{2n_i}},
\end{equation*}
together with (\ref{eq23}), we have
\begin{equation*}
I_j^{D_t^\epsilon}((d(t),0');\xi(t))=4d(t)^{2+\sum_{i\in I}2n_i/m_i}\prod_{i\in J}d_i^\ast(t)^{2n_i}I_j^{\Phi\circ\Sigma_t(D_t^\epsilon)}(0;(\Phi\circ\Sigma_t)'(d(t),0')\xi(t))
\end{equation*}
for $j=1,2$. Now we take sufficiently small $\epsilon,\delta>0$ arbitrary. From Lemma \ref{L33}, there exists $t_0(\epsilon,\delta)>0$ such that
\begin{equation*}
\mathscr{E}_\delta\times\ds\prod_{i\in J}B_{n_i}(0,1-\delta)\subset\Phi\circ\Sigma_t(D_t^\epsilon)
\end{equation*}
for all $t\in(0,t_0(\epsilon,\delta))$, which implies that
\begin{align*}
&I_j^{\Phi\circ\Sigma_t(D_t^\epsilon)}(0;(\Phi\circ\Sigma_t)'(d(t),0')\xi(t))\\
&\ge I_j^{\mathscr{E}_\delta\times\prod_{i\in J}B_{n_i}(0,1-\delta)}(0;(\Phi\circ\Sigma_t)'(d(t),0')\xi(t))\\
&=(1-\delta)^{2(1+\sum_{i\in I}n_i/m_i+\sum_{i\in J}n_i)}I_j^{\mathscr{E}\times\prod_{i\in J}B_{n_i}(0,1)}(0;(\Phi\circ\Sigma_t)'(d(t),0')\xi(t)).
\end{align*}
Here, we have used the monotonicity in (\ref{eq21}) for the first inequality and (\ref{eq23}) for the last equality. Hence, for every $t\in(0,t_0(\epsilon,\delta))$,
\begin{equation}\label{eq38}
\ds\frac{I_j^{\mathscr{E}\times\prod_{i\in J}B_{n_i}(0,1)}(0;(\Phi\circ\Sigma_t)'(d(t),0')\xi(t))}{I_j^{\Phi\circ\Sigma_t(D_t^\epsilon)}(0;(\Phi\circ\Sigma_t)'(d(t),0')\xi(t))}\le\ds\frac1{(1-\delta)^{2(1+\sum_{i\in I}n_i/m_i+\sum_{i\in J}n_i)}}.
\end{equation}
Since, from Lemma \ref{L33} again, 
\begin{equation*}
\Phi\circ\Sigma_t(D_t^\epsilon)\subset\mathscr{E}^\epsilon\times\ds\prod_{i\in J}B_{n_i}\left(0,\frac{d_i^\epsilon(t)}{d_i^\ast(t)}\right)
\end{equation*}
for all $t\in(0,t_0(\epsilon,\delta))$, we have
\begin{align*}
&I_j^{\Phi\circ\Sigma_t(D_t^\epsilon)}(0;(\Phi\circ\Sigma_t)'(d(t),0')\xi(t))\\
&\le I_j^{\mathscr{E}^\epsilon\times\prod_{i\in J}B_{n_i}(0,d_i^\epsilon(t)/d_i^\ast(t))}(0;(\Phi\circ\Sigma_t)'(d(t),0')\xi(t))\\
&=\ds\prod_{i\in I}\left(\ds\frac{c_i}{c_i-\epsilon}\right)^{2n_i/m_i}\!\!\!\times\ds\prod_{i\in J}\left(\ds\frac{d_i^\epsilon(t)}{d_i^\ast(t)}\right)^{2n_i}\!\!\times I_j^{\mathscr{E}\times\prod_{i\in J}B_{n_i}(0,1)}(0;(\Phi\circ\Sigma_t)'(d(t),0')\xi(t)),
\end{align*}
which implies that
\begin{equation}\label{eq39}
\ds\frac{I_j^{\mathscr{E}\times\prod_{i\in J}B_{n_i}(0,1)}(0;(\Phi\circ\Sigma_t)'(d(t),0')\xi(t))}{I_j^{\Phi\circ\Sigma_t(D_t^\epsilon)}(0;(\Phi\circ\Sigma_t)'(d(t),0')\xi(t))}\ge\ds\frac1{\ds\prod_{i\in I}\left(\ds\frac{c_i}{c_i-\epsilon}\right)^{2n_i/m_i}\!\!\!\times\ds\prod_{i\in J}\left(\ds\frac{d_i^\epsilon(t)}{d_i^\ast(t)}\right)^{2n_i}}.
\end{equation}
Now 
\begin{equation}\label{eq310}
\begin{split}
&\ds\frac{4d(t)^{2+\sum_{i\in I}2n_i/m_i}\prod_{i\in J}d_i^\ast(t)^{2n_i}I_j^{\mathscr{E}\times\prod_{i\in J}B_{n_i}(0,1)}(0;(\Phi\circ\Sigma_t)'(d(t),0')\xi(t))}{I_j^{D_t}((d(t),0');\xi(t))}\\
&=\ds\frac{I_j^{\mathscr{E}\times\prod_{i\in J}B_{n_i}(0,1)}(0;(\Phi\circ\Sigma_t)'(d(t),0')\xi(t))}{I_j^{\Phi\circ\Sigma_t(D_t^\epsilon)}(0;(\Phi\circ\Sigma_t)'(d(t),0')\xi(t))}\\
&\ \ \ \ \ \times\ds\frac{4d(t)^{2+\sum_{i\in I}2n_i/m_i}\prod_{i\in J}d_i^\ast(t)^{2n_i}I_j^{\Phi\circ\Sigma_t(D_t^\epsilon)}(0;(\Phi\circ\Sigma_t)'(d(t),0')\xi(t))}{I_j^{D_t}((d(t),0');\xi(t))}\\
&=\ds\frac{I_j^{\mathscr{E}\times\prod_{i\in J}B_{n_i}(0,1)}(0;(\Phi\circ\Sigma_t)'(d(t),0')\xi(t))}{I_j^{\Phi\circ\Sigma_t(D_t^\epsilon)}(0;(\Phi\circ\Sigma_t)'(d(t),0')\xi(t))}\times\ds\frac{I_j^{D_t^\epsilon}((d(t),0');\xi(t))}{I_j^{D_t}((d(t),0');\xi(t))},
\end{split}
\end{equation}
where the last equality is due to (\ref{eq23}). Using (\ref{eq38}), (\ref{eq310}), and Lemma \ref{L34}, it follows that
\begin{equation*}
\begin{split}
\ds\limsup_{t\to+0}&\ds\frac{4d(t)^{2+\sum_{i\in I}2n_i/m_i}\prod_{i\in J}d_i^\ast(t)^{2n_i}I_j^{\mathscr{E}\times\prod_{i\in J}B_{n_i}(0,1)}(0;(\Phi\circ\Sigma_t)'(d(t),0')\xi(t))}{I_j^{D_t}((d(t),0');\xi(t))}\\
&\le\ds\frac1{(1-\delta)^{2(1+\sum_{i\in I}n_i/m_i+\sum_{i\in J}n_i)}}
\end{split}
\end{equation*}
for all sufficiently small $\delta>0$, which implies
\begin{equation}\label{eq311}
\begin{split}
\ds\limsup_{t\to+0}\ds\frac{4d(t)^{2+\sum_{i\in I}2n_i/m_i}\prod_{i\in J}d_i^\ast(t)^{2n_i}I_j^{\mathscr{E}\times\prod_{i\in J}B_{n_i}(0,1)}(0;(\Phi\circ\Sigma_t)'(d(t),0')\xi(t))}{I_j^{D_t}((d(t),0');\xi(t))}\le1.
\end{split}
\end{equation}
Also, using (\ref{eq39}), (\ref{eq310}), Lemma \ref{L32} and \ref{L34}, it follows that
\begin{equation*}
\begin{split}
\ds\liminf_{t\to+0}&\ds\frac{4d(t)^{2+\sum_{i\in I}2n_i/m_i}\prod_{i\in J}d_i^\ast(t)^{2n_i}I_j^{\mathscr{E}\times\prod_{i\in J}B_{n_i}(0,1)}(0;(\Phi\circ\Sigma_t)'(d(t),0')\xi(t))}{I_j^{D_t}((d(t),0');\xi(t))}\\
&\ge\ds\frac1{\ds\prod_{i\in I}\left(\ds\frac{c_i}{c_i-\epsilon}\right)^{2n_i/m_i}\times(1+\epsilon)^{\sum_{i\in J}4n_i/\mu_i}}
\end{split}
\end{equation*}
for all sufficiently small $\epsilon>0$, and hence,
\begin{equation}\label{eq312}
\begin{split}
\ds\liminf_{t\to+0}\ds\frac{4d(t)^{2+\sum_{i\in I}2n_i/m_i}\prod_{i\in J}d_i^\ast(t)^{2n_i}I_j^{\mathscr{E}\times\prod_{i\in J}B_{n_i}(0,1)}(0;(\Phi\circ\Sigma_t)'(d(t),0')\xi(t))}{I_j^{D_t}((d(t),0');\xi(t))}\ge1.
\end{split}
\end{equation}
By (\ref{eq311}) and (\ref{eq312}), we get
\begin{equation*}
\ds\lim_{t\to+0}\ds\frac{4d(t)^{2+\sum_{i\in I}2n_i/m_i}\prod_{i\in J}d_i^\ast(t)^{2n_i}I_j^{\mathscr{E}\times\prod_{i\in J}B_{n_i}(0,1)}(0;(\Phi\circ\Sigma_t)'(d(t),0')\xi(t))}{I_j^{D_t}((d(t),0');\xi(t))}=1.
\end{equation*}
A similar argument using the equation (\ref{eq22}) yields
\begin{equation*}
\ds\lim_{t\to+0}\ds\frac{4d(t)^{2+\sum_{i\in I}2n_i/m_i}\prod_{i\in J}d_i^\ast(t)^{2n_i}I_0^{\mathscr{E}\times\prod_{i\in J}B_{n_i}(0,1)}(0)}{I_0^{D_t}(d(t),0')}=1,
\end{equation*}
and therefore we obtain the desired statements.
\end{proof}

\section{Proof of the theorem}
We now give the proof of Theorem \ref{MT} by combining the lemmas proved above.

\begin{proof}[Proof of Theorem \ref{MT}]
Let $\xi\in\C^{n+1}\setminus\{0\}$. Recall the notations $\xi=\xi_{N,t}+\sum_{i=1}^k\xi_{T,t}^i$, $\xi_{T,t}^i\in T_{\pi(q(t))}^\C(bD)$ for $i=1,\dots,k$ and $\xi_{N,t}\bot T_{\pi(q(t))}^\C(bD)$.  As
\begin{equation*}
R_t^2\circ R_t^1(\xi)=\left\langle R_t^2\circ R_t^1(\xi),e_1\right\rangle e_1+\ds\sum_{i=1}^k\sum_{j=1}^{n_i}\left\langle R_t^2\circ R_t^1(\xi),e_j^i\right\rangle e_j^i,
\end{equation*}
we have
\begin{equation*}
\xi_{N,t}=\left\langle R_t^2\circ R_t^1(\xi),e_1\right\rangle(R_t^2\circ R_t^1)^{-1}(e_1),\text{ and }
\end{equation*}
\begin{equation*}
\xi_{T,t}^i=\ds\sum_{j=1}^{n_i}\left\langle R_t^2\circ R_t^1(\xi),e_j^i\right\rangle(R_t^2\circ R_t^1)^{-1}(e_j^i),\text{ for }i=1,\dots,k,
\end{equation*}
which implies 
\begin{equation}\label{eq41}
|\xi_{N,t}|=\left|\left\langle R_t^2\circ R_t^1(\xi),e_1\right\rangle\right|,\text{ and }
\end{equation}
\begin{equation}\label{eq42}
|\xi_{T,t}^i|=\left(\ds\sum_{j=1}^{n_i}\left|\left\langle R_t^2\circ R_t^1(\xi),e_j^i\right\rangle\right|^2\right)^{1/2},\text{ for }i=1,\dots,k.
\end{equation}
Since
\begin{equation*}
\cal{K}_D(q(t))=\cal{K}_{D_t}(d(t),0'),
\end{equation*}
\begin{equation*}
B_D(q(t);\xi)=B_{D_t}((d(t),0');R_t^2\circ R_t^1(\xi)),\text{ and }
\end{equation*}
\begin{equation*}
H_D(q(t);\xi)=H_{D_t}((d(t),0');R_t^2\circ R_t^1(\xi)),
\end{equation*}
it follows, from (\ref{eq24}) and Lemma \ref{L35}, that
\begin{equation}\label{eq43}
\ds\lim_{t\to+0}\ds\frac{4d(t)^{2+\sum_{i\in I}2n_i/m_i}\prod_{i\in J}d_i^\ast(t)^{2n_i}\cal{K}_D(q(t))}{\cal{K}_{\mathscr{E}\times\prod_{i\in J}B_{n_i}(0,1)}(0)}=1,
\end{equation}
\begin{equation}\label{eq44}
\ds\lim_{t\to+0}\ds\frac{B_D(q(t);\xi)}{B_{\mathscr{E}\times\prod_{i\in J}B_{n_i}(0,1)}(0;(\Phi\circ\Sigma_t)'(d(t),0')R_t^2\circ R_t^1(\xi))}=1,\text{ and }
\end{equation}
\begin{equation}\label{eq45}
\ds\lim_{t\to+0}\ds\frac{H_D(q(t);\xi)-2}{H_{\mathscr{E}\times\prod_{i\in J}B_{n_i}(0,1)}(0;(\Phi\circ\Sigma_t)'(d(t),0')R_t^2\circ R_t^1(\xi))-2}=1.
\end{equation}
Here, 
\begin{equation}\label{eq46}
\begin{split}
(\Phi\circ\Sigma_t)'(d(t),0')R_t^2\circ R_t^1(\xi)=-&\left\langle R_t^2\circ R_t^1(\xi),e_1\right\rangle\ds\frac{e_1}{2d(t)}\\
&+\ds\sum_{i\in I}\ds\sum_{j=1}^{n_i}\left\langle R_t^2\circ R_t^1(\xi),e_j^i\right\rangle\ds\frac{e_j^i}{d(t)^{1/m_i}}\\
&+\ds\sum_{i\in J}\ds\sum_{j=1}^{n_i}\left\langle R_t^2\circ R_t^1(\xi),e_j^i\right\rangle\ds\frac{e_j^i}{d_i^\ast(t)}.
\end{split}
\end{equation}
By using Lemma \ref{L26}, combining (\ref{eq41}), (\ref{eq42}), (\ref{eq43}), (\ref{eq44}), (\ref{eq45}) and (\ref{eq46}), the proof of the main theorem is complete.
\end{proof}

\section*{Acknowledgment}
The author would like to thank my advisor, Professor Joe Kamimoto, for carefully reading an earlier draft and for many helpful suggestions regarding the presentation of this paper.

\end{document}